\documentclass[11pt]{article}
\usepackage{amssymb}
\usepackage{url}
\newcommand{\be}{\begin{equation}}
\newcommand{\ee}{\end{equation}}
\newtheorem{lemma}{Lemma}
\begin{document}
\title{\textbf{The discrete logarithm problem over prime fields: the safe prime case. The Smart attack, non-canonical lifts and logarithmic derivatives}}
\author{H. Gopalakrishna Gadiyar and R. Padma \\
Department of Mathematics\\
School of Advanced Sciences\\ V.I.T. University, Vellore 632014 INDIA\\E-mail: gadiyar@vit.ac.in, rpadma@vit.ac.in}
\date{~~}
\maketitle

\begin{abstract}
In this brief note we connect the discrete logarithm problem over prime fields in the safe prime case to the logarithmic derivative.


\end{abstract}

\section{Introduction and the Main Idea} Let $a_0$ be a primitive root of  a prime number $p>2$. We know that for every $b_0 \in \left\{1,2,\cdots ,p-1\right\}$ there exists a unique integer $n_p$ modulo $p-1$ satisfying 
\be
a_0^{n_p} \equiv b_0 {\rm{~mod~}}p \, .\label{eq:basiceqn}
\ee
$n_p$ is called the discrete logarithm or index of $b_0$ to the base $a_0$ modulo $p$.
  In \cite{gadiyar:maini:padma} the authors got the Teichm{\"u}ller expansion using Hensel lifting of the discrete logarithm problem (\ref{eq:basiceqn}). This is got by raising both sides to the power $p$:
\be
a_0^{n_pp} \equiv b_0^p {\rm{~mod~}}p^2 
\ee
which can be written as
\be
(a_0+a_1p)^{n_p} \equiv (b_0+b_1p) {\rm{~mod~}}p^2 \, . \label{eq:teichp}
\ee
The Iwasawa logarithm of a $p$-adic number $x$ is defined as $\frac{1}{p-1} \log x^{p-1}$. 
As this vanishes for a Teichm{\"u}ller character the solution $n_p$ could not be found out, but a formula
\be
n_p \equiv \frac{(b_1-\beta_{n_p})/b_0}{a_1/a_0} {\rm{~mod~}}p
\ee
was obtained where $\beta_{n_p}$ is the carry 
\be
a_0^{n_p} \equiv b_0+\beta_{n_p} p {\rm{~mod~}}p^2 \, .
\ee
Kontsevich \cite{kontsevich} and Riesel \cite{riesel} point out that the difficulty arises because the problem is stated modulo $p$ and the solution is needed modulo $p-1$. Hence we go to the discrete logarithm problem modulo the composite modulus $p(p-1)$. In this connection, see Bach \cite{bach}.

In this paper we consider primes $p$ of the form $2q+1$ where $q$ is a prime number. $p$ is called a safe prime as it is believed that the discrete logarithm problem is computationally difficult in this case when $p$ is `large'. 

From (\ref{eq:basiceqn}) we can go to the discrete logarithm problem 
\be
a_0^n \equiv b_0 {\rm{~mod~}}pq \, .
\ee
(See Lemma 1.) From the assumptions made in Lemma 1 $a_0$ generates a subgroup of order $q\phi(q)$ modulo $pq$. Hensel lifting the problem modulo $p^2q^2$ we get
\be
(a_0+a_1pq)^n \equiv (b_0+b_1pq) {\rm{~mod~}}p^2 q^2\, . \label{eq:teichpq}
\ee
The order of the group generated by $a_0+a_1pq$ remains as $q\phi(q)$ modulo $p^2q^2$. Also
\be
(a_0+a_1pq)^{q\phi(q)} \equiv 1 {\rm{~mod~}}p^2 q^3\, . \label{eq:teichpower}
\ee
(See Lemma 2.)  
Expanding (\ref{eq:teichpq}) using the binomial theorem, we get
\be
a_0^n + n a_0^{n-1} a_1 pq \equiv b_0 + b_1 pq {\rm{~mod~}}p^2 q^2\, .
\ee
Writing
\be
a_0^n \equiv b_0 + \beta_n pq {\rm{~mod~}}p^2 q^2\, ,
\ee
will give
\be
\beta_n +n \frac{b_0}{a_0} a_1 \equiv b_1 {\rm{~mod~}}pq \, . \label{eq:formulapq}
\ee
Here $\beta_n$ is the carry of $a_0^n$ modulo $p^2q^2$ and note that $n$ and $\beta_n $ are the two unknowns in the above linear congruence.

The summary of what we have done so far is that there are three problems when we try to solve the discrete logarithm problem modulo $p$: 
\begin{enumerate}
\item[1.] The problem is given modulo $p$ and the solution is needed modulo $p-1$.
\item[2.] The Iwasawa logarithm of the Teichm{\"u}ller expansion modulo $p^2$ is 0
\item[3.]The binomial theorem on the Teichm{\"u}ller expansion modulo $p^2$ gives 'carry'.
\end{enumerate}

We overcome the first problem by going modulo $pq$. The fact that we cannot get $n$ arises from two possibilities being blocked as in the modulo $p$ case. The analogue of the Teichm{\"u}ller expansion does not have a non-zero logarithm (see (\ref{eq:teichpower})) and if the binomial theorem is used, a carry occurs as in the case of mod $p$, see (\ref{eq:formulapq}). 

However if we can construct a non-canonical lift modulo $p^2q^2$ then the problems dissolve. Thus solving the discrete logarithm problem is equivalent to the construction of a non-canonical lift. 

The non-canonical lifts exist and can be written in the form 
\be
(a_0+(a_1+k)pq)^n \equiv (b_0+(b_1+l))pq) {\rm{~mod~}}p^2 q^2\, . \label{eq:noncan}
\ee
When $k=k_1p$ for some $k_1 \not\equiv 0$ mod $q$, then $l=l_1p$ for some $l_1$mod $q$. In this case the order of the group is $q\phi(q)$. For the other $k$ and $l$ modulo $pq$ the order of the group will be $pq\phi(q)$. On expanding (\ref{eq:noncan}) using the binomial theorem, one gets
\be
(a_0+a_1pq)^n~+~n(a_0+a_1pq)^{n-1}~ k~pq \equiv (b_0+b_1pq)~+~l~pq {\rm{~mod~}}p^2 q^2
\ee
and using (\ref{eq:teichpq})
\be
n \equiv \frac{l_1/b_0}{k_1/a_0} {\rm{~mod~}}q \, .
\ee
in the first case and
\be
n \equiv \frac{l/b_0}{k/a_0} {\rm{~mod~}}pq \, .
\ee
in the second case. 

If we use the notation $da_0$ for $k_1$ and $db_0$ for $l_1$ then 
\be
n \equiv \frac{db_0/b_0}{da_0/a_0} {\rm{~mod~}}q \, .  \label{eq:nmodq}
\ee
and if we use the notation $da_0$ for $k$ and $db_0$ for $l$ then 
\be
n \equiv \frac{db_0/b_0}{da_0/a_0} {\rm{~mod~}}pq \, . \label{eq:nmodpq}
\ee
Thus $n$ can be thought of as the logarithmic derivative. The non-canonical extensions (modulo $p^2q^2$) of the subgroup generated by $a_0$ mod $pq$ are labeled by $da_0$. As $p=2q+1$, once we get $n$ mod $q$, $n$ mod $p-1$ would be either $n$ or $n+q$ mod $p-1$.

Note that we can get (\ref{eq:nmodq}) and (\ref{eq:nmodpq}) by raising (\ref{eq:noncan}) to the powers $q\phi(q)$ and $pq\phi(q)$ respectively. In the second case we get
\be
\left(\left(a_0+(a_1+k)pq\right)^{pq\phi(q)}\right)^n\equiv \left(b_0+(b_1+l))pq\right)^{pq\phi(q)} {\rm{~mod~}}p^3 q^3\, ,
\ee
which on expanding and using the notation in Section 2 will give
\be
1+n(q(a_0)+\frac{(a_1+k)}{a_0}\phi(q))p^2q^2) \equiv 1+(q(b_0)+\frac{(b_1+l)}{b_0}\phi(q))p^2q^2  {\rm{~mod~}}p^3 q^3\,.
\ee
Using the formula for $a_1$ and $b_1$ one gets (\ref{eq:nmodpq}). This way of getting $n$ is analogous to the attack on anomalous elliptic curves by Smart \cite{smart}, Semaev \cite{semaev}, Satoh and Araki \cite{satoh:araki}. 

We would like to comment that derivatives of numbers have been studied historically for a long time starting from  Kummer \cite{hilbert}, \cite{washington}, A. Weil (expanded by Kawada) \cite{kawada} and more recently by A. Buium \cite{buium}. Hence the problem which is standing in isolation studied only by cryptologists gets connected to mainstream algebra and number theory. This was a complete surprise to the authors which is why we have written this brief note to bring it to the attention of experts in these areas.

\section{Lemmas} We need some definitions and notations before we prove our lemmas.
In \cite{lerch} Lerch defined the Fermat quotient for a composite modulus. Let $x$ be such that $gcd(x,n)=1.$ Then $q(x)$ defined by
\be
x^{\phi(n)} \equiv 1+q(x)n {\rm{~mod~}}n^2 \, .
\ee
is called the Fermat quotient of $x$ modulo $n$. We do not use the Euler's $\phi $-function but we use Carmichael's $\lambda $ function. $\lambda (n)$ is defined as follows \cite{Cameron:Preece}. $\lambda(2)=1$, $\lambda(4)=2$ and
\be
\lambda(n)=~\left \{ \begin{array}{ll} \phi(p^r),~&{\rm ~if ~} n=p^r\\
                                          2^{r-2},&{\rm ~if ~} n=2^r, ~r\ge 3\\
																				lcm(\lambda(p_1^{r_1}), \lambda(p_2^{r_2}),\cdots , \lambda(p_k^{r_k})),& {\rm ~if ~} n=p_1^{r_1}p_2^{r_2}\cdots p_k^{r_k} \end{array} \right .
\end{equation}
When $n=p^2q^2$ where $p=2q+1$ $q$ is a prime, $\phi(p^2q^2)=2pq^2\phi(q)$ and $\lambda(p^2q^2)=pq\phi(q)$. In other words the order of the group of units modulo $p^2q^2$ is $\phi(p^2q^2)$ whereas the order of the largest cyclic group modulo $p^2q^2$ is $\lambda(p^2q^2)$. Hence we define $q(x)$ by the congruence
\be
x^{pq\phi(q)} \equiv 1+q(x)p^2q^2 {\rm{~mod~}}p^3q^3 \, . \label{eq:fermatpq}
\ee

\begin{lemma} Let $a_0$ be a primitive root of $p$ and $q$. Let $gcd(b_0,~q)=1$. Then the congruence $a_0^n \equiv b_0 {\rm{~mod~}}p$ can be extended to
\be
a_0^{n_p} \equiv b_0 {\rm{~mod~}}pq \, . \label{eq:basiceqn2}
\ee
if and only if the Legendre symbols
\be
\left(\frac{b_0}{p}\right)=\left(\frac{b_0}{q}\right). \label{eq:legendresym}
\ee 
\end{lemma}

\noindent{{\bf Proof.}} $a_0^n \equiv b_0 {\rm{~mod~}}pq $ if and only if
$$
\begin{array}{lccccl}
a_0^n &\equiv& a_0^{n_p} &\equiv &b_0 {\rm{~mod~}}p &{\rm{~and}} \\
a_0^n &\equiv& a_0^{n_q} &\equiv &b_0 {\rm{~mod~}}q &~~ \, .
\end{array}
$$
This happens if and only if
\begin{eqnarray}
n &\equiv & n_p {\rm{~mod~}}p-1 {\rm{~and}} \\
n &\equiv & n_q {\rm{~mod~}}q-1 \, .
\end{eqnarray}
This is possible if and only if
\be
2=gcd(p-1, q-1) |(n_p-n_q).
\ee
by Chinese Remainder theorem.
That is 
\be
n_p \equiv n_q {\rm{~mod~}}2 \, .
\ee
In other words $b_0$ is a quadratic residue or nonresidue modulo $p$ and $q$ simultaneously. That is $\left(\frac{b_0}{p}\right)=\left(\frac{b_0}{q}\right).$

\begin{lemma} If $a_0^n \equiv b_0 {\rm{~mod~}}pq $ holds then 
\be
(a_0+a_1pq)^n \equiv (b_0+b_1 pq) {\rm{~mod~}}p^2 q^2\, , \label{eq:EQQ}
\ee
where 
\be
a_1=-\frac{q(a_0))a_0}{\phi(q)} {\rm {~and~}} b_1=-\frac{q(b_0))b_0}{\phi(q)} {\rm{~mod~}}pq\, . \label{eq:Eqn}
\ee
\end{lemma}
\noindent{\bf Proof} We want $a_1$ and $b_1$ to satisfy (\ref{eq:Eqn}). Using the carry notation 
\be
a_0^{n} \equiv b_0+\beta_n pq {\rm{~mod~}}p^2q^2 \, , \label{eq:neweQ}
\ee
we get the equation
\be
\beta_n +n \frac{b_0}{a_0} a_1 \equiv b_1 {\rm{~mod~}}pq \, . \label{eq:neweq}
\ee
Taking the power $pq\phi(q)$ on both sides of (\ref{eq:neweQ})
\be
a_0^{npq\phi(q)} \equiv (b_0+\beta_n pq)^{pq\phi(q)} {\rm{~mod~}}p^3q^3
\ee
and using (\ref{eq:fermatpq}) we get
\be
nq(a_0) \equiv q(b_0)+\frac{\beta_n}{b_0}\phi(q) {\rm{~mod~}} pq \,. \label{eq:neweq1}
\ee
Comparing (\ref{eq:neweq}) and (\ref{eq:neweq1}) will give the desired values of $a_1$ and $b_1$.

\noindent{\bf Remark 1.} Note that $a_1$ and $b_1$ can be calculated in polynomial time and the order of $(a_0+a_1pq)$ is $q\phi(q)$ modulo $p^2q^2$.

\noindent{\bf Remark 2.} Note that the Legendre symbols in (\ref{eq:legendresym}) can be calculated in polynomial time.

\noindent{\bf Remark 3.} We are given $b_0$ mod $p$. If (\ref{eq:legendresym}) fails for the given $b_0$ we can check the same for $b_0+kp$ for $k=1,2,3 \cdots $ until the condition is satisfied or we can multiply $b_0$ by $a_0^k$ for some $k$ and check the condition. In the first case $n_p$ does not change and in the second case $n_p$ becomes $n_p+k$ modulo $p-1$ or

\noindent{\bf Remark 4.} We can take $b_0^2$ mod $pq$ and consider the new discrete logarithm problem 
\be
a_0^n \equiv b_0^2 {\rm{~mod~}}pq \, ,
\ee
or 

\noindent{\bf Remark 5.} We can even relax the conditions in Lemma 1 as in our earlier preprint \cite{gadiyar:padma} as follows. Let $gcd(a_0,q)=1$ and $gcd(b_0,q)=1$. Let $a_0$ be a primitive root of $p$ and let $a_0$ and $b_0$ satisfy $a_0^n \equiv b_0 {\rm{~mod~}}p$. Then
\be
a_0^{n\phi(q)} \equiv b_0^{\phi(q)} {\rm{~mod~}}pq \, .
\ee

In this case the formulae corresponding to (\ref{eq:nmodq}) and (\ref{eq:nmodpq}) would be
\be
n \equiv \frac{db_0}{da_0} {\rm{~mod~}}q \, .  
\ee
and
\be
n \equiv \frac{db_0/b_0^{\phi(q)}}{da_0/a_0^{\phi(q)}} {\rm{~mod~}}pq \, .  
\ee

\section{Conclusion} For the composites $p^2q^2$ the Euler function $\phi (p^2q^2)=2pq^2\phi(q)$ and the Carmichael function $\lambda (p^2q^2)=pq\phi(q) $  are not equal. Also  $\lambda (p^2q^2) |\phi (p^2q^2)$ and hence many non-canonical lifts exist. As is well known this would involve a suitable choice of polynomial for lifting. Recall that the polynomials are $x^{p-1}-1$ and $x^{pq\phi(q)}-1$ in the cases of Teichm{\"u}ller lifting modulo $p^2$ and $p^2q^2$ respectively. This attack can be generalized to elliptic curve discrete logarithm problem over prime fields where $q$ will be connected to the order of the group. See \cite{silverman} for various ways of lifting the elliptic curve discrete logarithm problem.

\end{document}